\input amstex
\documentstyle{amsppt}
\magnification\magstep1
\NoBlackBoxes
\def\a{\alpha}
\def\d{\delta}

\def\R{{\bold R}}
\def\o{\omega}

\define\cl#1{\overline {#1}}

\def\sB{{\Cal B}}

\def\sbs{\subset}

\def\e{\epsilon}

\def\ti{\times}

\def\supp{\operatorname{supp}}

\def\Card{\operatorname{Card}}

\newcount\secno
\secno=-1
\newcount\subsecno
\subsecno=-1
\newcount\thmno
\thmno=0
\newcount\bigthmno
\bigthmno=0
\newcount\secthmno
\bigthmno=0

\def\ifundefined#1{\expandafter\ifx\csname#1\endcsname\relax}

\def\privateName#1{\ifundefined{privateName{#1}}\errmessage{Undefined name
#1}\else\csname privateName{#1}\endcsname\fi}

\def\name#1{\privateName{#1}}

\def\resetName#1#2{{\global\expandafter\edef\csname
privateName{#1}\endcsname{#2}}}

\def\setName#1#2{\ifundefined{privateName{#1}}\resetName{#1}{#2}\else
\errmessage{Name #1 is already in use}\fi}

\def\theoremName#1{%
\global\advance\thmno by1%
\global\advance\bigthmno by1%
\global\advance\secthmno by1%
\setName{#1}{%
        \ifnum\number\the\secno<0\else
               \ifnum\number\the\subsecno<0\else\number\the\subsecno.\fi
               \number\the\secno.%
        \fi
        \number\the\thmno
        }%
}

\def\theoremLabel#1{setName{#1}{#1}}

\def\No#1{\theoremName{#1}\name{#1}}

\def\label#1{\theoremLabel{#1}\privateName{#1}}

\def\newSection{\ifnum\secno<0%
                        \secno=1%
                \else
                        \advance\secno by1%
                \fi
                \ifnum\subsecno<0\else
                \subsecno=1%
                \fi
                \thmno=0%
                \secthmno=0}

\def\newSubSection{\ifnum\subsecno<0\subsecno=1\else
\advance\subsecno by1\fi
\thmno=0}

\def\setTheoremNumber#1{\thmno=#1\advance\thmno by-1}

\def\showTotalTheoremNumber{\message{{\number\the\bigthmno} claims etc.
total}}

\def\showSectionTheoremNumber{\message{The section contains
{\number\the\bigthmno} claims etc.}}


\def\setBibl#1#2{\setName{*bibl*#1}{#2}}
\def\bibl#1{\name{*bibl*#1}}

\setBibl{Praha?}{7} \setBibl{Pol?}{6}
\setBibl{DU}{2} 
\setBibl{OP}{5}
\setBibl{AG}{1}
\setBibl{E}{3}
\setBibl{vM}{4}
\setBibl{P}{6}

\topmatter
\author  V. V. Uspenskij
\endauthor

\title A note on a question of R.~Pol concerning light maps
\endtitle

\abstract Let $f:X\to Y$ be an onto map between compact spaces such that
all point-inverses of $f$ are
zero-dimensional. Let $A$ be the set of all functions $u:X\to I=[0,1]$
such that $u[f^\leftarrow(y)]$ is zero-dimensional for all $y\in Y$.
Do almost all maps $u:X\to I$, in the sense of Baire category, belong to $A$?
H.~Toru\'nczyk proved that the answer is
yes if $Y$ is countable-dimensional. We extend this result to the case when
$Y$ has property $C$. 
\endabstract

\address
321 Morton Hall, Department of Mathematics, Ohio
University, Athens, Ohio 45701, USA
\endaddress

\email uspensk\@bing.math.ohiou.edu, vvu\@uspensky.ras.ru
\endemail




\endtopmatter

\document
In R.~Pol's article [\bibl{Pol?}]
the following question was posed [\bibl{OP}, Problem 423]:

\smallskip{\it
Let $f:X\to Y$ be a continuous map of a compactum $X$ onto a compactum $Y$
with $\dim f^\leftarrow(y)=0$ for all $y\in Y$. Does there exist a nontrivial
continuous function $u:X\to I$ into the unit interval such that
$u[f^\leftarrow(y)]$ is zero-dimensional for all $y\in Y$?
}\smallskip

It was shown in [\bibl{DU}] 
that the answer is positive. However, R.~Pol informed
me that he actually meant another question: do almost all maps
$u:X\to I$, in the sense of Baire category, have the property considered
above? H.~Toru\'nczyk gave a positive answer under the assumption that
$Y$ is countable-dimensional. The aim of the present note is 
to extend this result to the case when
$Y$ has property $C$. In the general case the question remains open.

A space $Y$ is a {\it $C$-space\/},
or has {\it property~C\/}, if for any sequence $\{\a_n:n\in \o\}$ of open
covers of $Y$ there exists a sequence $\{\mu_n:n\in\o\}$
of disjoint families of open sets in $Y$
such that each $\mu_n$ refines $\a_n$ and
the union  $\bigcup_{n\in\o} \mu_n$ is a cover of $Y$.
This notion was first defined about~1973
by W.E.~Haver for compact metric spaces
and then by D.F.~Addis and J.H.~Gresham [\bibl{AG}] in the general case.
Every finite-dimensional
paracompact space and every countable-dimensional metric space has
property~$C$ [\bibl{AG}], [\bibl{E}].
Every normal $C$-space is weakly infinite-dimensional [\bibl{AG}], [\bibl{E}]
(Engelking [\bibl{E}] includes normality in the definition of property~$C$).
R.~Pol's example of a weakly infinite-dimensional compact
metric space which is not countable-dimensional ([\bibl{E}], [\bibl{vM}])
has property~$C$ and thus
distinguishes between property~$C$ and the property of being
countable-dimensional.
It is an open problem whether for (compact) metric
spaces property~$C$ is equivalent to the property of being
weakly infinite-dimensional.

All maps are assumed to be continuous.
We denote by $I$ the interval $[0,1]$.
A map $f:X\to Y$ is {\it light\/} if the fibres $f^\leftarrow(y)$ are
zero-dimensional for all $y\in Y$. For a compact space $X$ let $C(X,I)$ be
the space of all maps
$u:X\to I$ with the usual metric, induced by the metric of the 
Banach space $C(X)=C(X,\R)$.

\proclaim{\No{1}. Theorem}
Let $f:X\to Y$ be an onto light map between compact spaces. 
Let $A$ be the set of all functions $u:X\to I$
such that $u[f^\leftarrow(y)]$ is zero-dimensional for all $y\in Y$.
If $Y$ has property~$C$, then $A$ is a dense $G_\d$-subset of $C(X,I)$.
\endproclaim

The proof is based on a selection theorem for $C$-spaces obtained in
[\bibl{Praha?}].
Let $X$ be a paracompact $C$-space. Suppose that to each $x\in X$ 
a contractible non-empty subset $\Phi(x)$ of a space $Y$ is assigned. Suppose
that the multi-valued map 
$\Phi$ has the following semi-continuity property: for every compact
$K\sbs Y$ the set $\{x\in X: K\sbs \Phi(x)\}$ is open. Then $\Phi$ has
a continuous selection: there exists a continuous map $\phi: X\to Y$ such that 
$\phi(x)\in \Phi(x)$ for each $x\in X$ [\bibl{Praha?}, Theorem~1.3].
We shall use
a corollary of this theorem involving
the notion of a $Z$-set. 
Denote by $C(X,Y)$ the space of all maps $f:X\to Y$ in the
compact-open topology. 
Let us say that a closed subset $F$
of a topological space $X$ is a {\it $Z$-set\/} in $X$ if for any compact
space $K$ the set $C(K, X\setminus F)$ is dense in $C(K, X)$. If
$X$ is a separable metric ANR, this definition agrees with the usual
one~[\bibl{vM}]. 
For a closed subset
$F\sbs X$ to be a $Z$-set, it suffices that the identity map
of $X$ be in the closure of the subspace $C(X,X\setminus F)$ of $C(X,X)$.
If $F$ is a $Z$-set in $X$ and $U$ is open in $X$, then $F\cap U$
is $Z$-set in $U$.
If $C$ is a convex subset of a Banach space and
$F$ is a $Z$-set in $C$, then $C\setminus F$ is contractible (see, for example,
Proposition~6.4 in [\bibl{Praha?}]).
Therefore, the selection theorem formulated
above implies

\proclaim{\No{Th2}. Theorem}
Let X be a paracompact $C$-space. Let $C$ be a convex subset of a Banach space.
Suppose that to each $x\in X$ a $Z$-subset $Z(x)$ of $C$ and a convex subset
$U(x)$ of $C$ are assigned so that
the set $\bigcup_{x\in X} \{x\}\ti Z(x)$ is closed in $X\ti C$ and the set 
$\bigcup_{x\in X} \{x\}\ti U(x)$ is open in $X\ti C$. Then there
exists a continuous map $f:X\to C$ such that $f(x)\in U(x)\setminus Z(x)$ 
for every $x\in X$.
\endproclaim

It follows from the arguments of [\bibl{Praha?}] that Theorem~\name{Th2}
actually characterizes $C$-spaces among paracompact spaces. 

\proclaim{\No{zset}. Lemma}
Let $X$ be a convex subset  of a locally convex space $E$. If $Y$ is a convex
dense subspace of $X$, then any closed $F\sbs X$ which is disjoint from $Y$
is a $Z$-set in $X$.
\endproclaim

\demo{Proof}
It suffices to prove that for any convex symmetric
neighbourhood $V$ of zero in $E$
and for any compact $K\sbs X$ there exists a map $f:X\to Y$ such that 
$f(x)\in x+V$ for every $x\in K$. 
Since $Y$ is dense in $X$, we have $Y+V\supset X\supset K$, and by 
the compactness of $K$ there exists a finite $A\sbs Y$ such that $K\sbs A+V$.
Let $\{h_a:a\in A\}$ be a partition of unity subordinated to the cover
$\{a+V:a\in A\}$ of $K$. This means that each $h_a$ is a map from $K$ to $I$,
$\sum_{a\in A}h_a=1$ and the support $\supp(h_a)$ of $h_a$ is contained 
in $a+V$ for every $a\in A$.
The partition
of unity $\{h_a\}$ defines a map of $K$ into a simplex of dimension
$\Card(A)-1$.
This map can be extended over $X$, since a simplex is an absolute retract. 
It follows that there exists a partition
of unity $\{H_a:a\in A\}$ on $X$ such that the restriction of $H_a$ to
$K$ coincides with $h_a$ for every $a\in A$.
Define $f:X\to E$ by $f(x)=\sum_{a\in A}H_a(x)a$. The range of $f$ is 
contained in the convex hull of $A$ and hence in $Y$. Let us show that 
$f(x)-x\in V$ for every $x\in K$. Fix $x\in K$, and let $B$ be the set of all
$a\in A$ such that $h_a(x)>0$. If $a\in B$, then $x\in \supp(h_a)\sbs a+V$.
Therefore $f(x)-x=\sum_{\a\in A}h_a(x)(a-x)=\sum_{\a\in B}h_a(x)(a-x)$ is a 
convex combination of points of $V$ and hence belongs to $V$. 
\qed\enddemo

\proclaim{\No{zerodim}. Lemma}
Let $X$ be compact, $Y$ be a zero-dimensional closed subspace of $X$. If
$F$ is a closed subspace of $C(X,I)$ such that $f(Y)$ is infinite for
every $f\in F$, then $F$ is a $Z$-set in $C(X,I)$.
\endproclaim

\demo{Proof}
In virtue of Lemma~\name{zset}, it suffices to show that the convex set
$\{g\in C(X,I): g(Y)\text{ is finite}\}$ is dense in $C(X,I)$. Fix
$f\in C(X,I)$ and $\e>0$. Since $\dim Y=0$, there exists a map $h:Y\to I$
with finite range such that $0\le f(y)-h(y)<\e$ for every $y\in Y$. Let
$k:X\to I$ be an extension of the map $y\mapsto f(y)-h(y)$ 
$(y\in Y)$ over $X$ such that $k(x)<\e$ for every $x\in X$. The function 
$g=f-k\in C(X)$ is $\e$-close to $f$ and coincides with $h$ on $Y$, hence
$g(Y)=h(Y)$ is finite. If the range of $g$ is not contained in $I$, 
replace $g$ by $rg$, where $r$ is the natural retraction of the real line onto
$I$.
\qed\enddemo

\demo{Proof of Theorem~\name{1}}
Let $f:X\to Y$ be a light map of a compact space $X$ onto a compact $C$-space
$Y$. Let $A$ be the set of all maps $u:X\to I$ 
such that $u[f^\leftarrow(y)]$ is zero-dimensional for all $y\in Y$.
Let $C=C(X,I)$. 
We must show that $A$ is a dense $G_\d$-subset of $C$. 

For every subset $V\sbs I$ let
$A_V$ be the set of all maps $u:X\to I$ such that for every $y\in Y$
the set $u[f^\leftarrow(y)]$ does not contain $V$.
Fix a countable base $\sB$ in $I$. Since a subset of $I$ is zero-dimensional
if and only if it does not contain any element of $\sB$, we have 
$A=\bigcap_{V\in \sB} A_V$. Thus it suffices to prove that for every 
$V\in \sB$ the set $A_V$ is open and dense in $C$.

We show
that for every $V\sbs I$ the set $A_V$ is open in $C$. For every
$t\in I$ let $B_t$ be the set of all pairs $(y,u)$ in $Y\ti C$
such that $t\in u[f^\leftarrow(y)]$, and let $C_t$ be the set of all triples
$(x,y,u)$ in $X\ti Y\ti C$ such that $f(x)=y$ and $u(x)=t$.
Every $B_t$ is closed, since 
$B_t$ is the image of the closed set $C_t$ under
the projection $X\ti Y\ti C\to Y\ti C$
which is a closed map. Similarly, the
projection $Y\ti C\to C$ is closed and sends the closed set
$\bigcap_{t\in V} B_t$ to the complement of $A_V$. Hence $A_V$ is open in
$C$.

We prove that $A_V$ is dense in $C$ for every infinite subset 
$V\sbs I$. Fix $h\in C$ and $\e>0$. We must show that there exists
$w\in A_V$ such that $|w(x)-h(x)|<\e$ for every $x\in X$. For every $y\in Y$
let $U(y)$ be the convex set of all $u\in C$ such that $|u(x)-h(x)|<\e$
for every $x\in f^\leftarrow(y)$, and let $Z(y)$ be the set of all 
$u\in C$ such that $V\sbs u[f^\leftarrow(y)]$.
According to Lemma~\name{zerodim}, $Z(y)$ is a $Z$-set in $C$. 
The set
$\bigcup_{y\in Y}\{y\}\ti Z(y)$ is closed in $Y\ti C$, since it is equal
to the closed set $\bigcap_{t\in V} B_t$ considered in the preceding paragraph.
The set $\bigcup_{y\in Y}\{y\}\ti U(y)$
is open in $Y\ti C$, since its complement
is equal to the image of the closed subset $\{(x,u):|u(x)-h(x)|\ge\e\}$ of
$X\ti C$ under the perfect map $f\ti \text{id}_C: X\ti C\to Y\ti C$. 
Thus we can apply Theorem~\name{Th2}. In virtue
of this theorem, there exists a continuous map $y\mapsto u_y$ from $Y$ to $C$
such that $u_y\in U(y)\setminus Z(y)$ for every $y\in Y$. The map $w:X\to I$
defined by $w(x)=u_{f(x)}(x)$ has the required properties: $w\in A_V$ and
$|w(x)-h(x)|<\e$ for every $x\in X$. Indeed, for every $y\in Y$ the map
$w$ coincides with $u_y$ on the set $f^\leftarrow(y)$. Since $u_y\notin Z(y)$,
it follows that $w[f^\leftarrow(y)]=u_y[f^\leftarrow(y)]$ does not contain $V$.
Thus $w\in A_V$. Similarly, for every $x\in X$ we have $|w(x)-h(x)|=
|u_{f(x)}(x)-h(x)|<\e$, since $u_{f(x)}\in U(f(x))$. 
\qed\enddemo

\enddocument
\bye